\documentclass[12pt]{amsart}

\usepackage{amsmath} \usepackage{amsthm} \usepackage{mathrsfs}
\usepackage{latexsym}
\usepackage{amssymb}
\usepackage{enumerate}

\newtheorem{thm}{Theorem}[section]
\newtheorem*{thm*}{Theorem}
\newtheorem{fact}[thm]{Fact}
\newtheorem{conjecture}[thm]{Conjecture}
\newtheorem{claim}[thm]{Claim}
\newtheorem{cor}[thm]{Corollary}
\newtheorem*{cor*}{Corollary}

\newcommand\Fcal{\mathscr{F}}
\newcommand\Hcal{\mathscr{H}}

\newcommand{\tr}{\operatorname{tr}}
\newcommand{\cf}{\operatorname{cf}}

\newcommand{\Ord}{\mathrm{Ord}}

\newcommand\axiom{\mathrm}
\newcommand\SCH{\axiom{SCH}}
\newcommand\PFA{\axiom{PFA}}
\newcommand\MM{\axiom{MM}}
\newcommand\MRP{\axiom{MRP}}
\newcommand\ZFC{\axiom{ZFC}}

\renewcommand{\>}{\rangle}

\renewcommand{\epsilon}{\varepsilon}

\keywords{minimal walks, MRP, PFA, Prikry forcing, reflection, SCH}
\subjclass[2000]{03E05, 03E10, 03E75}

\title[{P}{F}{A}, {P}rikry forcing, and {S}{C}{H}]
{The {P}roper {F}orcing {A}xiom, {P}rikry forcing,
and the {S}ingular {C}ardinals {H}ypothesis}

\author{Justin Tatch Moore}
\thanks{The research presented
in this paper was supported by NSF grant DMS-0401893.}
\begin{document}

\begin{abstract}
The purpose of this paper is to present some results which suggest that
the Singular Cardinals Hypothesis follows from the Proper Forcing Axiom.
What will be proved is that a form of simultaneous reflection follows from
the Set Mapping Reflection Principle, a consequence of $\PFA$.
While the results fall short of showing that $\MRP$ implies $\SCH$,
it will be shown that $\MRP$ implies that if $\SCH$ fails first at
$\kappa$ then every stationary subset of
$S_{\kappa^+}^\omega = \{\alpha< \kappa^+:\cf(\alpha) = \omega\}$
reflects.
It will also be demonstrated that $\MRP$ always fails
in a generic extension by Prikry forcing.
\end{abstract}
\maketitle

\section{Introduction}

A stationary subset $S$ of a regular cardinal $\theta$ is said to
\emph{reflect} if there is a $\delta < \theta$ of uncountable
cofinality such that $S \cap \delta$ is stationary in $\delta$.
Similarly, a family $\Fcal$ of stationary sets is said to \emph{simultaneously
reflect} if there is a $\delta < \theta$ of uncountable
cofinality such that
$S \cap \delta$ is stationary in $\delta$ for every $S$ in $\Fcal$.
Notice the cofinality of $\delta$ acts as an upper bound for the number
of disjoint subsets of $\theta$ which can simultaneously reflect at
$\delta$.

Reflection and simultaneous reflection have been widely
studied in set theory with a number of
applications to areas such as
cardinal arithmetic, descriptive set theory,
and infinitary combinatorics.
Our starting point will be
the following theorem of Foreman, Magidor,
and Shelah.
\begin{thm}\cite{FMS}\label{MM_SCH}
Martin's Maximum implies that for every uncountable regular cardinal
$\theta > \omega_1$ and every collection $\Fcal$ of $\omega_1$
many stationary subsets of $S_\theta^\omega =
\{\alpha < \theta:\cf(\alpha) = \omega\}$
there is a $\delta < \theta$ of cofinality $\omega_1$ which simultaneously
reflects every element of $\Fcal$.
Moreover, it can be arranged that the union of $\Fcal$ contains a club
in $\delta$.
\end{thm}
Since for every regular uncountable $\theta$
there is a partition of $S_\theta^\omega$ into
disjoint stationary sets,
they conclude that $\MM$ implies that $\theta^{\omega_1} = \theta$
for all regular $\theta \geq \omega_2$.
By Silver's theorem \cite{SCH:Silver}
this in turn implies the Singular Cardinals Hypothesis
--- that $2^{\kappa} = \kappa^+$ for every singular strong limit $\kappa$.

In this paper, we will introduce and
explore a new notion of reflection called
\emph{trace reflection} and prove a result analogous to
Theorem \ref{MM_SCH}.
\begin{thm}
($\MRP$) Suppose that $\Omega \subseteq S_\theta^\omega$ is a non-reflecting
stationary set and that $\vec{C}$ avoids $\Omega$.
If $\Fcal$ is a collection of stationary subsets of $\Omega$ and $\Fcal$
has size $\omega_1$ then there is a $\delta < \theta$ of cofinality $\omega_1$
such that every element of $\Fcal$ simultaneously trace reflects at $\delta$.
\end{thm}
It will follow that $\MRP$ implies any failure of 
$\SCH$ must occur first at a singular cardinal $\kappa$
such that every stationary subset of $S^\omega_{\kappa^+}$ reflects.
The above theorem also has the following corollary.

\begin{cor} \label{MRP_fails}
Suppose that $M \subseteq V$ is an inner model with the same
cardinals such that for some cardinal $\kappa$
\begin{enumerate}

\item $\cf(\kappa)^V = \omega < \cf(\kappa)^M = \kappa$ and

\item every stationary subset of $\kappa^+$ in $M$ is stationary in $V$.

\end{enumerate}
Then $\MRP$ fails in $V$.
In particular, $\MRP$ fails in any generic extension by Prikry
forcing.\footnote{When I submitted this paper I was under the impression
that it was unknown whether $\PFA$ always failed in a Prikry extension.
Since the acceptence of this paper I have been made aware that this was
not the case.
Magidor has shown in an unpublished note that, by a slight modification
of an argument of Todor\v{c}evi\'{c},
$\PFA$ implies $\square_{\kappa,\omega_1}$ fails for all
$\kappa > \omega_1$.
On the other hand, Cummings and Schimmerling have
shown in \cite{indexed_squares} that after Prikry forcing at $\kappa$,
$\square_{\kappa,\omega}$ and hence $\square_{\kappa,\omega_1}$
always holds.}
\end{cor}

This paper is intended to be self contained.
Section \ref{tracereflection} contains the definition of trace reflection
and all of the necessary background on Todor\v{c}evi\'{c}'s trace function.
Section \ref{MRP} provides the necessary background on the Set Mapping
Reflection Principle which will figure prominently in the analysis.
The main results then follow in Section \ref{mainresults}.

The notation in the paper is mostly standard.
All ordinals are von Neumann ordinals --- the set of their
predecessors.
$H(\theta)$ is the collection of all sets of hereditary cardinality
less than $\theta$.
If $X$ is an uncountable set, $[X]^\omega$ is used to denote all
countable subsets of $X$.
See \cite{set_theory:Jech} or
\cite{set_theory:Kunen} for more background;
see \cite{multiple_forcing} for some information on the combinatorics
of $[X]^\omega$, the club filter, and stationary subsets of
$[X]^\omega$.

I would like to thank the referee for their careful reading useful
comments and suggestions. 
\section{Trace Reflection}
\label{tracereflection}

In this section I will define trace reflection.
First recall Todor\v{c}evi\'{c}'s notion of a walk on a given
cardinal $\theta$ (see \cite{partitioning_ordinals} or \cite{cseq}).
A \emph{$C$-sequence} is a sequence
$\vec{C} = \<C_\alpha:\alpha < \theta\>$ where
$\theta$ is an ordinal, $C_\alpha$ is closed and cofinal in $\alpha$ for
limit ordinals $\alpha$, and $C_{\alpha+1} = \{\alpha\}$.
A $C$-sequence $\vec{C}$ is said to \emph{avoid} a subset
$\Omega \subseteq \theta$ if $C_\alpha$ is disjoint from $\Omega$ for every
limit $\alpha < \theta$.
Notice that if $\Omega \subseteq \theta$ is a non-reflecting stationary set
then there is a $C$-sequence on $\theta$ which avoids $\Omega$.
Conversely, any $\Omega \subseteq \theta$ which is avoided by a $C$-sequence
cannot reflect.

For a given $C$-sequence, the \emph{trace function} is defined
recursively by
$$\tr(\alpha,\alpha) = \{\}$$
$$\tr(\alpha,\beta) = \tr(\alpha,\min(C_\beta \setminus \alpha)) \cup
\{\beta\}.$$
Hence $\tr(\alpha,\beta)$ contains all ordinals ``visited'' in the walk
from $\beta$ down to $\alpha$ along the $C$-sequence
except for the destination $\alpha$.\footnote{The omission of the
destination is not standard, but it simplifies the presentation at some
points.
For instance, Fact \ref{trace_fact} requires this omission.}
The following property of the trace function captures some of its most
important properties.
\begin{fact}\label{trace_fact}
If $\alpha < \beta$ and $\alpha$
is a limit then there is an $\alpha_0 < \alpha$ such that
$$\tr(\alpha,\beta) \subseteq \tr(\gamma,\beta)$$
whenever $\alpha_0 < \gamma < \alpha$.
If $\vec{C}$ avoids $\{\alpha\}$ then it can further be arranged that
$$\tr(\alpha,\beta) \cup\{\alpha\} \subseteq \tr(\gamma,\beta).$$
\end{fact}
\begin{proof}
First, observe that if $\xi$ is in $\tr(\alpha,\beta)$ then
either $C_\xi \cap \alpha$ bounded or else $\alpha$ is in
$C_\xi$.
Furthermore, the latter can only occur if $\xi$ is the least element
of $\tr(\alpha,\beta)$.
If $\alpha_0 < \alpha$ is an upper bound for every set
$C_\xi \cap \alpha$ such that $\xi$ is in $\tr(\alpha,\beta)$ and
$\alpha \not \in C_\xi$, then it is easily checked that $\alpha_0$
has the desired properties (use induction on $\beta$).
Such a bound exists since $\tr(\alpha,\beta)$ is finite.
Finally, if $\vec{C}$ avoids $\{\alpha\}$ then $\alpha$ is not
in $C_\xi$ for any $\xi$ in $\tr(\alpha,\beta)$.
It is therefore possible to prove the stronger conclusion in this case.
\end{proof}
Let $\theta$ be an ordinal of uncountable cofinality.
For a given $C$-sequence $\vec{C}$ of length $\theta$, define $\Hcal(\vec{C})$ to be
the collection of all $X \subseteq \theta$ such that whenever
$E \subseteq \theta$ is closed and unbounded,
there are $\alpha < \beta$ in $E$
with $\tr(\alpha,\beta) \cap X \ne \emptyset$.
Clearly the complement of $\Hcal(\vec{C})$ is a $\sigma$-ideal.

We say an element $X$ of $\Hcal(\vec{C})$ \emph{trace reflects
with respect to $\vec{C}$} if there is a $\delta < \theta$ of uncountable
cofinality such that $X \cap \delta$ is in
$\Hcal(\vec{C} \restriction \delta)$.
Simultaneous trace reflection is defined in a similar manner.
If $\vec{C}$ is clear from the context, I will omit the phrase
``with respect to $\vec{C}$.''
As with ordinary simultaneous reflection,
if $\Fcal$ is a disjoint family of elements of $\Hcal(\vec{C})$ which
simultaneously trace reflect at $\delta$, then the cardinality of
$\Fcal$ is at most the cofinality of $\delta$.

\section{Set Mapping Reflection}
\label{MRP}
Now I will recall some definitions associated with the
Set Mapping Reflection Principle.
For the moment let $X$ be a fixed uncountable set and let
$\theta$ be a regular cardinal such that $H(\theta)$ contains
$[X]^\omega$.
The set $[X]^\omega$ is equipped with a natural topology ---
the Ellentuck Topology ---
defined by declaring intervals of the form
$$[x,N] = \{Y\in [X]^\omega:x \subseteq Y \subseteq N\}$$
to be open where $x$ is a finite subset of $N$.
If $M$ is a countable elementary submodel of $H(\theta)$ and
$\Sigma$ is a subset of $[X]^\omega$ then we say that $\Sigma$
is \emph{$M$-stationary} if $E \cap \Sigma \cap M$ is non-empty whenever
$E \subseteq [X]^\omega$ is a closed unbounded set in $M$.
If $\Sigma$ is set mapping defined on a collection of countable
elementary submodels of $H(\theta)$ then we say that
$\Sigma$ is \emph{open stationary} if $\Sigma(M) \subseteq [X]^\omega$
is open and $M$-stationary for all relevant $M$.

A set mapping $\Sigma$ as above \emph{reflects} if
there is a continuous $\in$-chain $\<N_\nu:\nu < \omega_1\>$
in the domain of $\Sigma$ such that for every limit $\nu > 0$,
$N_\xi \cap X$ is in $\Sigma(N_\nu)$ for coboundedly many $\xi$
in $\nu$.
If this happens then $\<N_\nu:\nu < \omega_1\>$ is called a \emph{reflecting
sequence} for $\Sigma$.
The axiom $\MRP$ asserts that every open stationary set mapping
defined on a club reflects.
In \cite{MRP} it is shown that $\MRP$ is a consequence of $\PFA$.
It is also shown there
that it implies $2^{\aleph_0} = 2^{\aleph_1} = \aleph_2$ and that
$\square(\kappa)$ fails for all regular $\kappa> \omega_1$.

\section{The main results}
\label{mainresults}
We now proceed to the proof of the main theorem.

\begin{thm*} ($\MRP$)
Suppose that $\Omega \subseteq S_\theta^\omega$ is a non-reflecting
stationary set and that $\vec{C}$ avoids $\Omega$.
If $\Fcal$ is a collection of stationary subsets of $\Omega$ and $\Fcal$
has size $\omega_1$ then there is a $\delta < \theta$ of cofinality $\omega_1$
such that every element of $\Fcal$ simultaneously trace reflects at $\delta$.
\end{thm*}

\begin{proof}
Let $\Fcal = \{\Omega_\xi: \xi < \omega_1\}$ be given and let
$\{S_\xi:\xi < \omega_1\}$ be a sequence of disjoint
stationary sets such that $\xi < \min(S_\xi)$ and
$\bigcup_{\xi < \omega_1} S_\xi$ contains a club.
For $M$ a countable elementary submodel of $H({2^\theta}^+)$ which contains
$\Fcal$, define
$\Sigma_\Fcal(M)$ to be the collection of all countable $N \subseteq \theta$
such that either $N \cap \theta$ has a last element
or else $\sup N < \sup (M \cap \theta)$ and
$$\tr(\sup N,\sup (M \cap \theta)) \cap \Omega_\delta \ne \emptyset$$
where $\delta$ is such that $M \cap \omega_1$ is in $S_\delta$.
That $\Sigma_{\Fcal}(M)$ is open is a consequence of Fact \ref{trace_fact}.

\begin{claim}
$\Sigma_{\Fcal}(M)$ is $M$-stationary.
\end{claim}

\begin{proof}
Let $E \in M$ be a club of countable subsets of $\theta$ and let
$\delta$ be such that $M \cap \omega_1$ is in $S_\delta$.
By elementarity and assumption that $\Omega_\delta$ is stationary, there
is an $\alpha$ in $\Omega_\delta \cap M$ such that for every $\alpha_0 < \alpha$,
there is an $N$ in $E \cap M$ such that $\alpha_0 < \sup(N) < \alpha$.
By Fact \ref{trace_fact} it is possible to find an $N$ in $E \cap M$ 
such that $\alpha$ is in $\tr(\sup(N), \sup(M \cap \theta))$.
\end{proof}

Now, let $\<N_\xi:\xi < \omega_1\>$ be a reflecting sequence for
$\Sigma_{\Fcal}$ and put
$$E = \{\sup(N_\xi \cap \theta): \xi < \omega_1\}$$
$$\delta = \sup E.$$
It suffices to show that for every $\xi < \omega_1$ and closed unbounded
$E' \subseteq \delta$ that
there are $\alpha < \beta$ in $E'$ such that
$\tr(\alpha,\beta) \cap \Omega_\xi$ is non-empty.
Let $\beta < \delta$ be a limit point of
$E \cap E'$ such that $N_\nu \cap \omega_1$ is in $S_\xi$ where
$\nu < \omega_1$ is such that $\beta = \sup (N_\nu \cap \theta)$.
By virtue of $\<N_\xi:\xi < \omega_1\>$ reflecting $\Sigma_{\Fcal}$
and
the definition of $E$, there is a $\beta_0 < \beta$ such that
if $\alpha$ is in $E$ with $\beta_0 < \alpha < \beta$ then
$\tr(\alpha,\beta) \cap \Omega_\xi$ is non-empty.
Selecting $\alpha$ in $E \cap E'$ with $\beta_0 < \alpha < \beta$,
we now have $\alpha < \beta$ both in $E'$ with
$\tr(\alpha,\beta) \cap \Omega_\xi$ non-empty as desired.
\end{proof}

We finish the section with proof of the corollary.
\begin{cor*}
Suppose that $M \subseteq V$ is an inner model with the same
cardinals such that for some cardinal $\kappa$
\begin{enumerate}

\item $\cf(\kappa)^V = \omega < \cf(\kappa)^M = \kappa$ and

\item every stationary subset of $\kappa^+$ in $M$ is stationary in $V$.

\end{enumerate}
Then $\MRP$ fails in $V$.
In particular, $\MRP$ fails in any generic extension by Prikry forcing.
\end{cor*}

\begin{proof}
Let $\vec{C}$ be a $C$-sequence in $M$ of length $\theta = \kappa^+$
such that for every $\alpha < \theta$, $C_\alpha$ has ordertype at most
$\kappa$ and $\vec{C}$ avoids
$$\Omega = \{\alpha < \theta:\cf(\alpha)^M = \kappa\}.$$
Let $\{\Omega_\xi:\xi < \kappa\}$ be a partition in $M$ of $\Omega$ into
disjoint stationary sets.
Pick an $X \subseteq \kappa$ in $V$ which is countable and
cofinal in $\kappa$.
Now suppose towards a contradiction that $\MRP$ holds in $V$.
By the main theorem there would be a $\delta < \theta$
of cofinality $\omega_1$ such that $\Omega_\xi$ trace reflects at $\delta$
for every $\xi$ in $X$.
Now observe that the cofinality of $\delta$ must be less than $\kappa$
in $M$ since otherwise it would have countable cofinality in $V$.
Let $E \subseteq \delta$ be closed and unbounded with
$|E| < \kappa$ and $E$ in $M$.
Put $$X^* = \{\xi < \kappa:
\exists \alpha, \beta \in E
(\tr(\alpha,\beta) \cap \Omega_\xi \ne \emptyset)\}.$$
Certainly $X^*$ is in $M$, has size less than $\kappa$ (since the
$\{\Omega_\xi:\xi < \kappa\}$ are all pairwise disjoint), and
is cofinal (since it contains $X$).
But this is a contradiction since $\kappa$ is regular in $M$.
\end{proof}

\section{Concluding remarks}

In Corollary \ref{MRP_fails}, it seems unlikely that the condition
on preserving cardinals or stationary sets is really necessary.
\begin{conjecture}
If $M$ is an inner model of $V$ 
such that
\begin{enumerate}

\item $2^{\omega_1} \cap M = 2^{\omega_1} \cap V$ and

\item $\Ord^{\omega_1} \cap M \ne \Ord^{\omega_1} \cap V$,

\end{enumerate}
then $\MRP$ fails in $V$.
\end{conjecture}

This seems very closely related to the next conjecture.

\begin{conjecture}
$\MRP$ implies that $\theta^{\omega_1} = \theta$ for all regular
$\theta \geq \omega_2$.
\end{conjecture}
\noindent
This would of course show that $\SCH$ follows from $\MRP$.

It should be remarked that Assaf Sharon has recently announced that
$\SCH$ can fail at $\kappa$ (even for $\kappa = \aleph_\omega$) and
yet every stationary subset of $\kappa^+$ reflects.
Hence it is not possible to prove the conjecture by establishing
a $\ZFC$ connection between the existence of a non-reflecting stationary
subset of $\kappa^+$ and the failure of $\SCH$ at $\kappa$. 
A possible approach, however, is to try to replace the assumption of a
non-reflecting stationary subset of $S^\omega_{\kappa^+}$
with the existence of a good scale for $\kappa$.
The motivating factor is that
a good scale for $\kappa$ always exists if $\SCH$
fails first at $\kappa$ (see Main Claim
1.3, p 46 in \cite{cardinal_arithmetic}).
One can also attempt to refute the existence of good scales
using $\MRP$ and thus prove $\MRP$ implies $\SCH$
(Magidor has shown that $\MM$ implies
good scales do not exist).
These approaches were suggested by
Veli\v{c}kovi\'{c} and Kojman respectively.

Finally, there are some results which link reflection in $[\lambda]^\omega$
to $\SCH$.
In \cite{FA_stationary} Veli\v{c}kovi\'c showed that if
$\theta > \omega_1$ is regular and stationary subsets of
$[\theta]^\omega$ reflect to an internally closed unbounded set
(strongly reflect in the language of \cite{FA_stationary}) then
$\theta^\omega = \theta$.
An immediate consequence is that $\PFA^+$ implies $\SCH$.
Recently Shelah improved this result by showing that
reflection of stationary subsets of $[\theta]^\omega$ to sets of
size $\omega_1$ is already sufficient to deduce $\theta^\omega = \theta$
\cite{reflection_SCH}.


\begin{thebibliography}{10}

\bibitem{indexed_squares}
James Cummings and Ernest Schimmerling.
\newblock Indexed squares.
\newblock {\em Israel J. Math.}, 131:61--99, 2002.

\bibitem{FMS}
Matthew Foreman, Menachem Magidor, and Saharon Shelah.
\newblock {M}artin's {M}aximum, saturated ideals, and nonregular ultrafilters.
 {I}.
\newblock {\em Ann. of Math. (2)}, 127(1):1--47, 1988.

\bibitem{multiple_forcing}
T.~Jech.
\newblock {\em Multiple forcing}, volume~88 of {\em Cambridge Tracts in
 Mathematics}.
\newblock Cambridge University Press, Cambridge, 1986.

\bibitem{set_theory:Jech}
Thomas Jech.
\newblock {\em Set theory}.
\newblock Perspectives in Mathematical Logic. Springer-Verlag, Berlin, second
 edition, 1997.

\bibitem{set_theory:Kunen}
Kenneth Kunen.
\newblock {\em An introduction to independence proofs}, volume 102 of {\em
 Studies in Logic and the Foundations of Mathematics}.
\newblock North-Holland, 1983.

\bibitem{MRP}
Justin~Tatch Moore.
\newblock Set mapping reflection.
\newblock submitted to JML in Nov. 2003.

\bibitem{cardinal_arithmetic}
Saharon Shelah.
\newblock {\em Cardinal arithmetic}, volume~29 of {\em Oxford Logic Guides}.
\newblock The Clarendon Press Oxford University Press, New York, 1994.
\newblock Oxford Science Publications.

\bibitem{reflection_SCH}
Saharon Shelah.
\newblock Reflection implies {S}{C}{H}.
\newblock Preprint, July 2004.
Sh:794. arXiv:math.LO/0404323.

\bibitem{SCH:Silver}
Jack Silver.
\newblock On the singular cardinals problem.
\newblock In {\em Proceedings of the International Congress of Mathematicians
 (Vancouver, B. C., 1974), Vol. 1}, pages 265--268. Canad. Math. Congress,
Montreal, Que., 1975.

\bibitem{partitioning_ordinals}
Stevo Todor\v{c}evi\'{c}.
\newblock Partitioning pairs of countable ordinals.
\newblock {\em Acta Math.}, 159(3--4):261--294, 1987.

\bibitem{cseq}
Stevo Todor\v{c}evi\'{c}.
\newblock Coherent sequences.
\newblock In {\em Handbook of Set Theory}. North-Holland, (in preparation).

\bibitem{FA_stationary}
Boban Veli{\v{c}}kovi{\'c}.
\newblock Forcing axioms and stationary sets.
\newblock {\em Adv. Math.}, 94(2):256--284, 1992.

\end{thebibliography}
\end{document}